\documentclass[12pt]{amsart}
\usepackage{graphicx,amssymb}
\usepackage[all]{xy}
\newtheorem{thm}{Theorem}[section]
\newtheorem{cor}[thm]{Corollary}
\newtheorem{lem}[thm]{Lemma}
\newtheorem{prop}[thm]{Proposition}
\theoremstyle{definition}

\newtheorem{ex}[thm]{Example}
\newtheorem{rem}[thm]{Remark}
\newtheorem{hyp}[thm]{Hypothesis}
\numberwithin{equation}{section}

\newcommand{\cK}{\mathcal{K}}

\newcommand{\cW}{\mathcal{W}}
\newcommand{\cV}{\mathcal{V}}
\newcommand{\cD}{\mathcal{D}}

\newcommand{\QQ}{\mathbb Q}

\newcommand{\ra}{\rightarrow}

\DeclareMathOperator{\Ima}{{Im}}

\DeclareMathOperator{\Core}{{Core}}
\DeclareMathOperator{\Gal}{{Gal}}

\DeclareMathOperator{\End}{{End}}

\def\Z{{\mathbb Z}}
\def\Q{{\mathbb Q}}

\def\P{{\mathbb P}}
\def\C{{\mathbb C}}

\def\cU{{\mathcal U}}
\def\G{{\mathcal G}} 
\def\bs{{\backslash}}

\DeclareMathOperator{\tr}{tr}
\DeclareMathOperator{\im}{Im}

\begin{document}

\title[ ]{Prym-Tyurin varieties using self-products of groups}%
\author{A. Carocca, H. Lange, R. E. Rodr{\'\i}guez and A. M. Rojas}

\address{A. Carocca\\Facultad de Matem\'aticas, Pontificia Universidad Cat\'olica de Chile, Casilla 306-22, Santiago, Chile}

\address{R. E. Rodr{\'\i}guez\\Facultad de Matem\'aticas, Pontificia Universidad Cat\'olica de Chile, Casilla 306-22, Santiago, Chile}

\address{H. Lange\\Mathematisches Institut,
              Universit\"at Erlangen-N\"urnberg\\Germany}

\address{A. M. Rojas \\Departamento de Matem\'aticas, Facultad de
Ciencias, Universidad de Chile, Santiago\\Chile}
\email{ }%

\thanks{The first, third and fourth author were supported by Fondecyt
grants 1060743, 1060742 and 11060468 respectively}%
\subjclass{14H40, 14K10,  }%
\keywords{Prym-Tyurin variety, correspondence}%

\begin{abstract}
Given Prym-Tyurin varieties of exponent $q$ with respect to a
finite group $G$, a subgroup $H$ and a set of rational irreducible
representations of $G$ satisfying some additional properties, we
construct a Prym-Tyurin variety of exponent $[G:H]q$ in a natural
way. We study an example of this result, starting from the
dihedral group $\mathbf{D}_p$ for any odd prime $p$. This
generalizes the construction of \cite{LRR} for $p=3$. Finally, we
compute the isogeny decomposition of the Jacobian of the curve
underlying the above mentioned example.
\end{abstract}
\maketitle

\section{Introduction}

Consider a cartesian diagram of {\it smooth} projective curves
\begin{equation} \label{diag1.1}
\xymatrix{
      &  X   \ar[dl]^{}_{q_1} \ar[dr]_{}^{q_2}\\
     X_1  \ar[dr]_{h_1}^{}  & & X_2  \ar[dl]^{h_2}\\
    &  \mathbb{P}^1    }
\end{equation}
with $\deg h_i = d_i \geq 2$. Then the Jacobian $JX$ cannot be equal to $JX_1 \times JX_2$, since the genus $g_X$ of $X$ is always bigger than the sum $g_{X_1} + g_{X_2}$. However, if the $h_i$ do not factorize via a cyclic \'etale covering for $i =1$ and 2, then $q_i^*: JX_i \ra JX$ is injective for $i = 1$ and 2 and it is easy to see that
$$
q_1^* + q_2^*: JX_1 \times JX_2 \ra JX
$$
is an isogeny onto its image. One might ask whether $q_1^* + q_2^*$ is an embedding. This is certainly not the case if $d_1$ is not a divisor of $d_2$ (we assume that $d_1 \leq d_2$), since then the type of the polarization on $JX_1 \times JX_2$ induced by the canonical polarization would be $(d_1, \ldots, d_1,d_2, \ldots, d_2)$.
However, if we assume $d_1 = d_2 =: d$, then the polarization would be the $d$-fold of a principal polarization and $JX_1 \times JX_2$ would be a Prym-Tyurin variety of exponent $d$ in $JX$. Recall that a {\it Prym-Tyurin variety of exponent} $q$ in a Jacobian $J$ is by definition an abelian subvariety of $J$ on which the canonical polarization of $J$ induces the $q$-fold of a principal polarization.

An example for this was given by Mumford in \cite[p. 346]{M}, where he showed that in the case of the above diagram with $\deg h_i =2$ for $i = 1$ and 2, if $Y$ denotes the hyperelliptic curve ramified over the union of the branch loci of $h_1$ and $h_2$, then there is an \'etale double covering $X \ra Y$ whose Prym variety in $JX$ is isomorphic (as principally polarized abelian varieties) to the product of Jacobians $JX_1 \times JX_2$.

Our main result is a generalization of Mumford's theorem. In order to state it, we need to recall some facts of \cite{clrr}. In that paper we started with a smooth projective curve $Z$ with action by a finite group $G$ such that $Z/G \simeq \P^1$. Then we associated to a set of pairwise non-isomorphic irreducible rational representations and a subgroup $H$ of $G$ satisfying some additional properties, a Prym-Tyurin variety $P$ in $JX$ with $X = Z/H$.
We say that this construction of the Prym-Tyurin variety is a {\it presentation of} $P$ with respect to the action of the group $G$, of the subgroup $H$ and the set of representations. For a more precise statement see Theorem \ref{thm4.8} below and the subsequent definition. To state our result, we also need the following definition.

 Let $\pi: Z \ra Y$ be a Galois covering with Galois group $G$ and let $C_1, \ldots, C_t$ be pairwise different conjugacy classes of cyclic subgroups of $G$. We define the {\it geometric signature} of the covering $\pi$ to be the tuple $[\gamma,(C_1,s_1),\ldots, (C_t,s_t)]$, where $\gamma$ is the genus of the quotient curve $Y$, the covering has a total of $\sum_{j=1}^t s_j$ branch points in $Y$ and exactly $s_j$ of them are of type $C_j$ for $j = 1, \ldots, t$; that is, the corresponding points in its fiber have stabilizer belonging to $C_j$.
Using this terminology, our first result can roughly be stated as follows.
\begin{thm} \label{thm1.1}
Suppose that for $i =1$ and 2 the action of a group $G$ on a curve $Z_i$ has geometric signature $[0;(C_1^i,s_1^i), \ldots, (C_{t_i}^i,s_{t_i}^i)]$ and the branch loci of the corresponding coverings are disjoint in $\P^1$. Assume that we are given a presentation of a Prym-Tyurin variety $P_i$ (with the same exponent $q$) with respect to the action of $G$ on $Z_i$, a subgroup $H$ and a set of non-trivial rational representations. Then

{\em (a):} the induced action of the group $G^2 := G \times G$ on the fibre product $Z := Z_1 \times_{\P^1} Z_2$ defines a Prym-Tyurin variety
$P$ in the Jacobian $J(Z/H^2)$ of exponent $\widetilde{q} =[G:H]q$;

{\em (b):} under mild additional assumptions, \quad $P \simeq P_1 \times P_2$ \; (isomorphic as principally polarized abelian varieties).
\end{thm}
For a more precise version of this theorem see Theorems \ref{prod} and \ref{thm3.4}. In \cite{clrr} we showed that almost all known Prym-Tyurin varieties admit presentations in the above sense. Moreover, it is easy to find presentations of Prym-Tyurin varieties of exponent 1 (which are Jacobians according to the Theorem of Matsusaka-Ran) and exponent 2 (which essentially are classical Prym varieties according to a Theorem of Welters). In particular, this gives the above mentioned sufficient criterion for a product of Jacobians to be a subvariety of a Jacobian. Hence Theorem \ref{thm1.1} can be considered as a generalisation of Mumford's Theorem mentioned above. It should be noted that the special case of our proof in Mumford's situation is different from Mumford's original proof (see Example \ref{ex3.7} and Corollary \ref{cor3.9} below). \\

In the second part of the paper we work out an explicit example, which
was, in fact,  the starting point of this paper. In \cite{LRR} a Prym-Tyurin variety of exponent 3 was associated to every \'etale degree 3 covering $\widetilde{Y}$ of a hyperelliptic curve $Y$, such that the Galois group of its Galois closure over $\P^1$ is $\mathbf{S}_3 \times \mathbf{S}_3$, where $\mathbf{S}_3$ denotes the symmetric group of degree $3$. Solomon generalized in \cite{rudy} this construction to define a Prym-Tyurin variety for every group $\mathbf{S}_n \times \mathbf{S}_n$. Since $\mathbf{S}_3$ coincides with the dihedral group $\mathbf{D}_3$ of order 6, one might ask whether there is also a generalization of the construction of \cite{LRR} to $\mathbf{D}_n \times \mathbf{D}_n$. We show in the second part of the paper that this is in fact the case, at least for $\mathbf{D}_p$ with $p$ an odd prime. The construction certainly also generalizes to the group $\mathbf{D}_n$, however we restrict ourselves to a prime number, since the group theory of an arbitrary $\mathbf{D}_n$ is more complicated.

Let $p$ be an odd prime number, and consider an  \'etale $p$-fold covering $\widetilde{Y} \ra Y$ of a hyperelliptic curve $Y$ of genus $g$, such that the Galois closure of the composed map $\widetilde{Y} \ra Y \ra \P^1$ has Galois group $\mathbf{D}_p \times \mathbf{D}_p$. Then the fibre product
$$
X = \P^1 \times_{Y^{(2)}} \widetilde{Y}^{(2)}
$$
is a smooth projective curve. Here $Y^{(2)}$ and $\widetilde{Y}^{(2)}$ denote the second symmetric product of $Y$ and $\widetilde{Y}$ and
$\P^1 \ra Y^{(2)}$ is the canonical embedding of the $g^1_2$ of $Y$.
We define an effective symmetric fixed-point free $(p-1)^2$-correspondence $D$ on $X$ whose associated endomorphism $\gamma_D \in \End(JX)$ satisfies the equation
$$
\gamma_D^2 + (p-2)\gamma_D - (p-1) =0.
$$
This gives part (a) of the following theorem (see Corollary \ref{cor6.2}).
\begin{thm} \label{thm1.2}
{\em (a):} $D$ defines a Prym-Tyurin variety $P$ of exponent $p$ in the Jacobian $JX$. \\
{\em (b):} There exist smooth projective curves $X_1$ and $X_2$ whose fibre product over $\P^1$ is $X$, and such that
$$
P \simeq JX_1 \times JX_2,
$$
isomorphic as principally polarized abelian varieties.
\end{thm}

This can be proven directly, which was our first approach. Then we realized that this is a special case of Theorem \ref{thm1.1}. In fact, one can associate to the data the curves $X_1$ and $X_2$ in a canonical way, such that
$X = X_1 \times_{\P^1} X_2$. Moreover it turns out that the Kanev correspondence, associated to the correspondence of Theorem \ref{thm1.1}, coincides with the correspondence $D$ (see Proposition \ref{prop6.1}). So Theorem \ref{thm1.1} implies Theorem 1.2.\\

In the last part of the paper we compute the isogeny decomposition of the Jacobian of the curve $X$ of Theorem \ref{thm1.2}. Let $X$ be as above. The Galois group of the Galois closure $Z$ over $\P^1$ is the group $\mathbf{D}_p^2$. Then, if $X_1$ and $X_2$ denote the curves of Theorem \ref{thm1.2}. we show in Section \ref{S:desc},

\begin{thm}
There are \'etale $p$-fold coverings $\widetilde{Y_j} \ra Y, \; j = 1, \ldots, \frac{p-1}{2},$ of a hyperelliptic curve $Y$, all of them subcovers of $Z$ by explicitly given subgroups of $\mathbf{D}_p^2$, such that
$$
JX \sim JX_1 \times JX_2 \times P(\widetilde{Y}_1/Y) \times \cdots \times P(\widetilde{Y}_{(p-1)/2}/Y),
$$
where $\sim$ denotes isogeny and $P(\widetilde{Y_j}/Y)$ denotes the (generalized) Prym variety of the covering $\widetilde{Y_j} \ra Y$.
\end{thm}

Throughout the paper we work over the field of complex numbers.

\section{Preliminaries}

In this section we recall the main result of \cite{clrr}. Let $G$ be a finite group. In order to fix the notation, we recall some basic properties of representations of $G$. For any complex irreducible representation $V$ of $G$, we denote by $\chi_V$ its character, by $L_V$ its field of definition and by $K_V$ the subfield $K_V = \Q(\chi_V(g) \, | \, g \in G )$. Then $L_V$ and $K_V$ are finite abelian extensions of $\Q$; we denote by $m_V = [L_V : K_V]$ the Schur index of $V$. For any automorphism $\varphi$ of $L_V/\Q$ we denote by $V^{\varphi}$ the representation conjugate to $V$ by $\varphi$.

If $\cW$ is a rational irreducible representation of $G$, then there exists a complex irreducible representation $V$ of $G$, uniquely determined up to conjugacy in $\Gal(L_V/\Q)$, such that
$$
\cW \otimes_{\Q} \C \simeq \oplus_{\varphi \in \Gal(L_V/\Q)} V^{\varphi} \simeq m_V \oplus_{\tau \in \Gal(K_V/\Q)} V^{\tau}.
$$
We call $V$ a complex irreducible representation {\it associated to $\cW$}.

Let $\cW_1, \ldots, \cW_r$ denote nontrivial pairwise non-isomorphic
rational irreducible representations of the group $G$ with associated
complex irreducible representations $V_1, \ldots, V_r$. In our applications $r$ will be either $1$ or $2$. We make the
following hypothesis on the $\cW_i$ and a subgroup $H$ of $G$.

\begin{hyp} \label{hyp}
For all $k,l = 1 ,\ldots, r$ we assume
\begin{enumerate}
  \item[a)] $\dim V_k = \dim V_l =: n$,
  \item[b)] $K_{V_k} = K_{V_l} =: L$,
  \item[c)] $\dim V_k^H = 1$,
  \item[d)] $H$ is maximal with property c); that is, for every subgroup $N$ of $G$ with $H \subsetneqq N$
there  is an index $k$ such that $\dim V_k^N = 0$.
\end{enumerate}

\end{hyp}

Choose a set of representatives
$$
\{ g_{ij} \in G \;| \; i = 1, \ldots, d \; \mbox{and} \; j = 1, \ldots, n_i \}
$$
for both the left cosets and right cosets of $H$ in $G$, such that
$$
G = \sqcup_{i=1}^d Hg_{i1}H \quad \mbox{and} \quad Hg_{i1}H = \sqcup_{j=1}^{n_i} g_{ij}H = \sqcup_{j=1}^{n_i} Hg_{ij}
$$
are the decompositions of G into double cosets, and of the double
cosets into right and left cosets of $H$ in $G$. Moreover, we assume
$g_{11} = 1_G$.\\

Now let $Z$ be a (smooth projective) curve with $G$-action and $\pi_H : Z \to X := Z/H$. In \cite{clrr} we defined a correspondence on $X$, which is given by
\begin{equation} \label{eqnD}
{\cD}(x) = \sum_{i=1}^d b_{i} \sum_{j=1}^{n_i} \pi_H g_{ij}(z).\
\end{equation}
for all $x \in X$ and $z \in Z$ with $\pi_H(z) = x$, where
$$
b_i := \sum_{k=1}^r \sum_{h \in H} \tr_{L/\QQ}(\chi_{V_k}(hg_{i1}^{-1}))
$$
is an integer for $i = 1, \ldots, d$.

Let $\delta_{\cD}$ denote the endomorphism of the Jacobian
$JX$ associated to the correspondence ${\cD}$. We denote
by
$$
P_{\cD} := \im(\delta_{\cD})
$$
the image of the endomorphism $\delta_{\cD}$ in the
Jacobian $JX$ and call it the {\it (generalized) Prym variety}
associated to the correspondence ${\cD}$.

Setting
\begin{equation} \label{eq3.5}
b := \gcd \{ b_1 - b_i \; | \; 2 \leq i \leq d \},
\end{equation}
\cite[Theorem 4.8]{clrr} can be stated as follows:
\begin{thm} \label{thm4.8}
Let $\cW_1, \ldots, \cW_r$
denote nontrivial pairwise non-isomorphic rational irreducible representations of the group $G$
with associated complex irreducible representations $V_1, \ldots, V_r$ satisfying Hypothesis \ref{hyp} for a subgroup $H$ of $G$.
Suppose that the action of the finite group $G$ has geometric signature $[0;(C_1,s_1), \ldots , (C_{t},s_{t})]$
satisfying
\begin{equation}  \label{eqn2.2}
\sum_{j=1}^t s_j \left[ q[L:\QQ] \left(\sum_{k=1}^r  (\dim V_k - \dim V_k^{G_j})\right) - \left( [G:H] -|H \backslash G/G_j|\right) \right] = 0.
\end{equation}
where $G_j$ is a subgroup of $G$ of class $C_j$ and
$$q = \frac{|G|}{b \cdot n}.
$$
Then $P_{\cD}$ is a Prym-Tyurin
variety of exponent $q$ in $JX$.
\end{thm}

Furthermore, we showed in \cite[Section 4.4]{clrr},
\begin{equation}
  \label{eq:dim}
\dim P_{\cD} = [L:\mathbb{Q}] \sum_{i=1}^r [ - n +\frac{1}{2} \sum_{j=1}^t
s_j (\dim V_i - \dim V_i^{G_j}) ]
\end{equation}
and
\begin{equation}
  \label{eq:genus}
g_{X} = 1 -[G:H] + \frac{1}{2} \sum_{j=1}^t s_j ([G:H] - |H\backslash G / G_j|)
\end{equation}
In the sequel we will use the following definition:  We say that
the construction of the Prym-Tyurin variety $P = P_{{\cD}}$ of
Theorem \ref{thm4.8} is a {\it presentation of} $P$ with respect
to the action of the group $G$, the subgroup $H$ and the set of
representations $\{\cW_1, \ldots, \cW_r\}.$

\section{The construction}

In this section we show how to construct new Prym-Tyurin varieties out of given ones. As above let $G$ be a finite group.
For $i=1,2$, let $Z_i$ denote a smooth projective curve, on which $G$ acts with geometric signature $[0;(C_1,s^1_1), \ldots , (C_{t_1}, s^1_{t_1})]$ on $Z_1$ and $[0;(C_1',s^2_1), \ldots , (C_{t_2}', s^2_{t_2})]$ on $Z_2$. Let $h_i: Z_i \ra \P^1$ denote the corresponding covering maps. Then we consider the product group
$$
G^2 := G \times G.
$$
For any subgroup $H$ of $G$, we denote $H_1 = H \times \{1\}_{G}, \; H_2 = \{1\}_{G} \times H$ and write $C^1_j = C_j \times \{1\}_{G}$ and similarly
$C^2_j = \{1\}_{G} \times C_j'$. So the $C^{\ell}_j$ are considered as conjugacy classes of cyclic subgroups of $G^2$ for $\ell = 1$ and $2$.

\begin{lem}\label{lem:Z}
Suppose the branch loci of $h_i: Z_i \ra \P^1$ are disjoint in
$\P^1$. Then the fibre product
$$
Z := Z_1 \times_{\P^1} Z_2
$$
is a smooth projective curve, Galois over $\P^1$ with Galois group
$G^2$ and geometric signature $[0;(C^1_1,s^1_1), \ldots , (C^1_{t_1}, s^1_{t_1}),(C_1^2,s^2_1), \ldots , (C_{t_2}^2, s^2_{t_2})]$.
\end{lem}

\begin{proof}
This is elementary. The branch loci being disjoint implies that $Z$ is smooth and classical Galois theory implies that $Z$ is Galois over $\P^1$ with Galois group $G^2$. The last statement is clear from the definitions.
\end{proof}
Now let $H$ be a subgroup of $G$ and consider $X_1 := Z_1/H$ and $X_2 := Z_2/H$. Clearly the curve $X:= Z/(H \times H)$ is the fibre product of $X_1$ and $X_2$ over $\P^1$,
$$
X = Z/H^2 = X_1 \times_{\P^1} X_2
$$
and we have the following diagram, where the maps are the obvious ones.
\begin{equation} \label{diag3.1}
\xymatrix@R=0.5cm{
                 &      &  Z \ar[ddll]_{p_1} \ar[d]_{\pi} \ar[ddrr]^{p_2}     \\
                 &     &   X \ar[dl]_{q_1} \ar[dr]^{q_2} \\
Z_1 \ar[r]^{\pi_1} \ar[drr]_{h_1} & X_1 \ar[dr]^{\varphi_1} & & X_2 \ar[dl]_{\varphi_2} & Z_2 \ar[l]_{\pi_2} \ar[dll]^{h_2}\\ & & \P^1     }
\end{equation}
Here
$$
\deg h_1 = \deg h_2 = \deg p_1 = \deg p_2 = |G|,
$$
$$
\deg \varphi_1 = \deg \varphi_2 = \deg q_1 = \deg q_2 = [G:H],
$$
$$
\deg \pi = |H|^2.
$$

In the rest of this section we change slightly the notation for the geometric signature, for the sake of simplicity. Let $C_1, \ldots, C_t$ denote {\it all} conjugacy classes of nontrivial cyclic subgroups of $G$. Then in the tuple $[0;(C_1,s_1), \ldots, (C_t,s_t)]$ a number $s_j$ will be $0$ precisely
if there is no branch point of type $C_j$. Before constructing a Prym-Tyurin variety using the action of $G^2$ on $Z$, we need the following Lemma for which use the following notation.
For $i = 1$ and 2 let ${\cD}_i$ denote the correspondence on the curve $X_i$ defined in \eqref{eqnD} for the group $G$ with respect to the subgroup $H$ and the representations $\cW_1, \ldots , \cW_r$. Similarly let ${\cD}$ denote the correspondence  on the curve $X$ defined in \eqref{eqnD} for the Group $G^2$ with respect to the subgroup $H^2$ and the representations $\cW^1_i = \cW_i \otimes V_0$ and $\cW^2_i = V_0 \otimes \cW_i$ for $i = 1, \ldots , r$, where $V_0$ denotes the trivial representation of $G$. Then we have

\begin{lem} \label{lem3.6}
$$
{\cD}=|H|\cdot (q_1^* {\cD}_1+q_2^* {\cD}_2)
$$
\end{lem}

\begin{proof} For simplicity we assume $r=1$ and write $\cW = \cW_1$ etc. The proof for the general case is the same, only notationally more complicated.

Set $d=|H\bs G/H|$ and $\{g_{ij}:i=1,\ldots,d; \; j=1, \ldots ,n_i\}$ as in Section 2. Therefore
$|H^2\bs G^2/H^2|=d^2$ and
$\{(g_{ij},g_{kl}):i,k = 1,\ldots, d;\; j=1, \ldots ,n_i,\; l=1, \ldots , n_k\}$ are
representatives of both left and right cosets of $H^2$ in $G^2$.
According to \eqref{eqnD} we have
$$
{\cD}_{\nu}(x_{\nu}) = \sum_{i=1}^d a_{i} \sum_{j=1}^{n_i} \pi_H
g_{ij}(z_{\nu})
$$
for all $x_{\nu} \in X_{\nu}$, where $z_{\nu} \in Z_{\nu}$ is a
preimage of $x_{\nu}$, $1\leq {\nu}\leq 2$, and where
$$
a_i := \sum_{h \in H} \tr_{L/\QQ}(\chi_{V}(hg_{i1}^{-1}))
$$
is the same integer for both ${\cD}_{\nu}$ and $i = 1, \ldots, d$.

By definition we have $(q_{\nu}^*{\cD}_{\nu})(x_1,x_2)=q_{\nu}^{-1}{ \cD}_{\nu}q_{\nu}(x_1,x_2)=q_{\nu}^{-1}{\cD}_{\nu}(x_{\nu})$ for ${\nu}=1,2$. Therefore
\begin{equation} \label{eqn3.5}
(q_{\nu}^*{\cD}_{\nu})(x_1,x_2)=q_{\nu}^{-1}(\sum_{i=1}^da_{i} \sum_{j=1}^{n_i} \pi_H g_{ij}(z_{\nu}))=\sum_{i,k=1}^da_i\sum_{j=1}^{n_i}\sum_{l=1}^{n_k}(\pi_H g_{ij}(z_1),\pi_H g_{kl}(z_2)).
\end{equation}
On the other hand, according to \eqref{eqnD} we have for ${\cD}$,
\begin{equation} \label{eqn3.6}
{\cD}(x_1,x_2) = \sum_{i=1}^d\sum_{k=1}^d (b_{ik} + b_{ik}') \sum_{j=1}^{n_i}\sum_{l=1}^{n_k}(\pi_H g_{ij}(z_1),\pi_H g_{kl}(z_2))
\end{equation}
with $z_1$ and $z_2$ as above,
\begin{equation} \label{eqn3.8}
\begin{array}{rl}
b_{ik} &=  \sum_{(h_1,h_2) \in H^2} \tr_{L/\QQ}(\chi_{(V \otimes V_0)}(h_1g_{i1}^{-1},h_2g_{k1}^{-1}))\\
& =\sum_{(h_1,h_2) \in H^2} \tr_{L/\QQ}(\chi_{V}(h_1g_{i1}^{-1})\chi_{V_0}(h_2g_{k1}^{-1}))\\
& =\sum_{(h_1,h_2) \in H^2} \tr_{L/\QQ}(\chi_{V}(h_1g_{i1}^{-1}))\\
& =|H|\sum_{h \in H} \tr_{L/\QQ}(\chi_{V}(hg_{i1}^{-1}))\\
&=|H|\cdot a_i
\end{array}
\end{equation}
and similarly
\begin{equation} \label{eqn3.9}
b_{ik}' = \sum_{(h_1,h_2) \in H^2} \tr_{L/\QQ}(\chi_{(V_0 \otimes
V)}(h_1g_{i1}^{-1},h_2g_{k1}^{-1})) = |H| \cdot a_k \, .
\end{equation}
Inserting the last equations into \eqref{eqn3.6} and comparing with \eqref{eqn3.5} we deduce
$$
{\cD}(x_1,x_2)=|H|\cdot ({q_1^*\cD}_{1} + {q_2^* \cD}_{2})(x_1,x_2)
$$
\end{proof}

\begin{thm} \label{prod}

Let $\cW_1, \ldots , \cW_{r}$ denote nontrivial pairwise non-isomorphic rational irreducible representations of the group $G$ with associated complex irreducible representations $V_1, \ldots , V_{r}$ satisfying Hypothesis \ref{hyp} with respect to a subgroup $H$ of $G$.\\
Suppose that the action of $G$ on $Z_i$ has geometric signature $[0;(C_1,s^i_1), \ldots , (C_{t},s^i_{t})]$ on $Z_i$ for $i=1,2$, and satisfies
\begin{equation}  \label{eqn3.1}
\sum_{j=1}^t s^i_j \left[ q[L:\QQ]\sum_{k=1}^r  (\dim V_k - \dim V_k^{G_j})  -
\left( [G:H] -|H \backslash G/G_j|\right) \right] = 0.
\end{equation}
for $i = 1,2$, where $G_j$ is a subgroup of type $C_j$ and $q =
\dfrac{|G|}{b \cdot n}$. Furthermore assume that the branch loci
of $Z_i \ra \P^1$ are disjoint in $\P^1$.

Then the action of the group $G^2$ on the curve $Z$ defines a
Prym-Tyurin variety $P$ in the Jacobian $JX$ of exponent $
\tilde{q} = [G:H] \, q$ and dimension
$$
\dim P = [L:\mathbb{Q}] \sum_{i=1}^r [ - 2n +\frac{1}{2} \sum_{j=1}^{t}
(s^1_j +s^2_j) (\dim V_i - \dim V_i^{G_j})].
$$
\end{thm}
\begin{proof}
First note that Hypothesis \ref{hyp} is satisfied for the subgroup
$H^2$ of $G^2$ and the representations
$\cW_j \otimes V_0, V_0 \otimes  \cW_j$, with associated complex representations $V^1_j := V_j \otimes V_0$ and $V^2_j := V_0 \otimes V_j$, where $V_0$ denotes the trivial representation of $G$.

To see this, notice that
$$
\dim (V^i_j)^{H \times H}  = \langle V^i_j , \rho_{H \times H}^{G \times G} \rangle_{G \times G}  =
\langle V_j , \rho_H^G \rangle = 1
$$
for all $i$ and $j$. The maximality of $H \times H$ with respect to this property is a consequence of the fact that for every $\cW_j$ both representations $\cW_j \otimes V_0$ and $V_0 \otimes \cW_j$ occur in $\rho_{H \times H}^{G \times G}$.

First, we need to compute the exponent $\tilde{q}$ this data determines. In order to do that, we need to compute the greatest common divisor $\widetilde{b}$ of the differences between the first coefficient $\widetilde{b_1}$ of $\cD$ with the others. Using Equations \eqref{eqn3.6}, \eqref{eqn3.8} and \eqref{eqn3.9} of Lemma \ref{lem3.6} we obtain that
$$\widetilde{b_1}=2|H|a_1, $$
The rest of the coefficients of $\cD$ are of the following types
$$|H|(a_1+a_j) \;\text{ or }\; 2|H|a_j$$
therefore the differences are of the following types
$$|H|(a_1-a_j) \;\text{ or }\; 2|H|(a_1-a_j)$$
hence the corresponding greatest common divisor $\widetilde{b}$ is $|H|$ times
the corresponding one $b$ for one copy of $G$. The new exponent is then computed as

$$\widetilde{q}=\frac{|G^2|}{\widetilde{b}\cdot \dim V_j^1}=\frac{|G|\cdot |G|}{|H|b\cdot \dim V_j}=[G:H]q$$

Now, the assertion follows from Theorem \ref{thm4.8} as soon as we show that
\begin{equation} \label{eqn3.2}
\sum_{\ell=1}^2 \sum_{j=1}^t s_j^{\ell} \left[ \tilde{q}[L:\QQ]\sum_{i=1}^2\sum_{k=1}^r
(\dim V_k^{i} - \dim (V_k^{i})^{G_j^{\ell}})  -
([G^2:H^2] -|H^2 \backslash G^2/G_j^{\ell}|) \right] = 0.
\end{equation}

\noindent where $G_j^{\ell}$ is of class $C_j^{\ell}$ with $C_j^{\ell}$ as in Lemma \ref{lem:Z}.

To see this, observe that
\begin{equation} \label{eq:h}
  [G^2:H^2] =[G:H]^2 \, , \  \text{ and  } \
|H^2 \backslash G^2 /G_j^{\ell}| = [G:H]|H\backslash G /G_j|.
\end{equation}
Moreover, we have
\begin{equation} \label{eq:vg}
\dim (V_k^i)^{G_j^{\ell}}  = \langle V_k^i , \rho_{G_j^{\ell}}^{G^2}\rangle   =
\left\{
  \begin{array}{ll}
     \langle V_k, \rho_{G_j}^G \rangle = \dim V_k^{G_j}, & \hbox{if $\ell = i$;} \\
     \dim V_k^i = \dim   V_k, & \hbox{otherwise.}
  \end{array}
\right.
\end{equation}

It now follows from \eqref{eq:h} and \eqref{eq:vg} that  the left
hand side of \eqref{eqn3.2} equals $[G:H]$ times the sum with $i=1$ and $i=2$ of the left
hand side of \eqref{eqn3.1}. But these are zero by assumption, which implies the assertion. Finally, the computation of the dimension is a consequence of equation \eqref{eq:dim}.
\end{proof}
For our main applications we need only the following special case of the above theorem.
\begin{cor} \label{cor3.3}
Let $\cW$ be a nontrivial rational irreducible representation of $G$, with associated complex irreducible representation $V$, such that the subgroup $H$ of $G$ is maximal with the property
$$
\dim V^H = 1.
$$
Suppose that the action of $G$ on $Z_i$ has geometric signature $[0;(C_1,s^i_1), \ldots , (C_{t},s^i_t)]$ for $i =1,2$, and satisfies equation \eqref{eqn3.1} with $r=1$ for $i = 1$ and $2$. Furthermore assume that the branch loci of $Z_i \ra \P^1$ are disjoint in
$\P^1$.

Then the action of the group $G^2$ on the curve $Z$ defines a Prym-Tyurin variety $P$ in the Jacobian $JX$ of exponent
$ \tilde{q} = [G:H]q$ and dimension
$$
\dim P = [K_V:\mathbb{Q}] [ - 2\dim V +\frac{1}{2} \sum_{j=1}^t
(s^1_j + s^2_j) (\dim V - \dim V^{G_j}) ].
$$
\end{cor}

\begin{rem}
One could generalize Theorem \ref{prod} to an $n$-fold self-product $G^n$ of $G$, a subgroup $H$ of $G$ and $n$ Galois coverings $h_i:Z \ra \P^1$ with Galois group $G$ and pairwise disjoint branch loci in $\P^1$. We omit the details.
\end{rem}

The fact that the Prym-Tyurin variety $P$ is constructed via a product of groups suggests that it is a product of Prym-Tyurin varieties, and in fact, this is the case, as we will show in the next theorem. Let the notation be as in Theorem \ref{prod}. According to Theorem \ref{thm4.8} the action of the group $G$ and its subgroup $H$ induce a Prym-Tyurin variety $P_i$ of exponent $q$ in $JX_i$ for $i = 1$ and 2. Let $P$ denote the Prym-Tyurin variety in $JX$ of Theorem \ref{prod}. Then we have:

\begin{thm} \label{thm3.4}
Suppose the coverings $q_i:X \ra X_i$ do not factorize via an \'etale cyclic covering for $i=1,2$. Then we have
$$
P \simeq P_1 \times P_2.
$$
\end{thm}
\begin{proof}
The map $q_1^* + q_2^*: JX_1 \times JX_2 \ra JX$ is an isogeny onto its image. According to Lemma \ref{lem3.6} it maps $P_1 \times P_2$ into $P$. From equation \eqref{eq:dim} and Theorem \ref{prod} we get
$$
\dim (P_1 \times P_2) = \dim P.
$$
Hence $q_1^* + q_2^*$ induces an isogeny $P_1 \times P_2 \ra P$. Now, since the maps  $q_i:X \ra X_i$ do not factorize via a cyclic \'etale covering, the canonical polarization of $JX$ induces a polarization of the same type on $P$ and $P_1 \times P_2$, namely the $([G:H]\cdot q)$-fold of a principal polarization. This implies that $q_1^* + q_2^*: P_1 \times P_2 \ra P$ is an isomorphism.
\end{proof}

\begin{ex} \label{ex3.7}
As a first example consider the cyclic group $G = \Z_2$ of order 2. There is only one nontrivial rational irreducible representation, the alternating representation $\cW$. It certainly satisfies Hypothesis \ref{hyp} for the trivial subgroup $H = \{0\}$. For $i = 1$ and 2 let $X_i$ be a hyperelliptic curve of genus $g_{X_i}$. We assume that the hyperelliptic coverings $f_i:X_i \ra \P^1$ have disjoint branch loci in $\P^1$, so that we are in the situation of diagram \eqref{diag1.1}. The fibre product $X = X_1 \times_{\P^1} X_2$ is Galois over $\P^1$ with Galois group the Klein group $G^2$, so that in diagram \eqref{diag3.1} the curves $Z =X$ and $Z_i = X_i$ coincide.
Moreover the group $G$ acts on $X_i$ with geometric signature $[0;(G,2g_{X_i}+2)]$ and satisfies equation \eqref{eqn2.2} with $q = 1$. So the Prym-Tyurin variety of the action coincides with the Jacobian $JX_i$. Moreover the assumptions of Theorem \ref{prod} are fulfilled and we obtain a Prym-Tyurin variety $P$ of exponent 2 for the group $G^2$ in the Jacobian $JX$.

Let $K$ denote the third subgroup of $G^2$, i.e. $K = \{(0,0),(1,1)\}$ and denote
$$
Y := X/K.
$$
\begin{prop} \label{prop3.8}
{\em (a)} The curve $Y$ is hyperelliptic of genus $g_C = g_{X_1} + g_{X_2} + 1$ and the map $X \ra Y$ is an \'etale double covering.\\
{\em (b)}  The Prym-Tyurin variety $P$ of the action of $G^2$ coincides with the classical Prym variety of the \'etale covering $X \ra Y$.
\end{prop}

\begin{proof}
(a): Set $H_1=G\times \{0\}=\{(0,0),(1,0)\}$ and $H_2=\{0\}\times
G=\{(0,0),(0,1)\}$. Then $G^2$ acts on $X$ with geometric
signature $[0;(H_1,2g_{X_2}+2), (H_2,2g_{X_1}+2), (K,0)]$. From
this we conclude that the double covering $X \ra Y$ is \'etale. In
fact, if $x \in X$ is a branch point of the covering $X \ra \P^1$,
it cannot be one of $X \ra Y$, since $\text{Stab}_{G^2}(x) \cap K
= \{ (0,0) \}$. It follows that $Y \ra \P^1$ is a double covering
ramified in all $2(g_{X_1} + g_{X_2}) +4$ points, which means that
$Y$ is hyperelliptic of genus $g_{X_1}+ g_{X_2} + 1$.

(b): First note that the Prym variety of the covering $X\to Y$ is given
by the correspondence $1-\iota$, where $\iota$ is the involution
on $X$ with quotient $Y$. Hence $\iota$ is given by the element $(1,1)$ of $G^2$.

On the other hand, the Prym-Tyurin variety $P$ for the action of
$G^2$ is defined by the correspondence ${\cD}$ (see \eqref{eqnD})
constructed using the sum of the representations $\cW_1=\cW\otimes
V_0$ and $\cW_2=V_0 \otimes \cW$ of $G^2$, where $V_0$ and $\cW$
are the trivial and the alternating representation of $G$
respectively. We have then $H^2=\{(0,0)\}$, the trivial subgroup
of $G^2$, and the double coset representatives are just the
elements of $G^2$, i.e. $g_{11}=(0,0),g_{21}=(1,0),
g_{31}=(0,1),g_{41}=(1,1)$. They are also the left and right coset
representatives of $H^2$ in $G^2$.

By definition, we have for $i=1, \ldots ,4$,
$$
b_i := \sum_{k=1}^2\sum_{h \in H^2} (\chi_{\cW_k}(hg_{i1}^{-1}))=\chi_{\cW_1}(g_{i1})+
\chi_{\cW_2}(g_{i1}),
$$
which gives $b_1=2,\;b_2=b_3=0$ and $b_4=-2$. Hence
$b := \gcd \{ b_1 - b_i \; | \; 2 \leq i \leq 4 \}=2$.
Therefore
$${\cD}=1\cdot (0,0)- 1\cdot (1,1).$$
Hence  the correspondences coincide, which implies the assertion.
\end{proof}

As an immediate consequence of Theorem \ref{thm3.4} and Proposition \ref{prop3.8} we obtain:

\begin{cor} \label{cor3.9}
{\em (Mumford)}\\
Let $X_i$ be a hyperelliptic curve of genus $g_i$ for $i=1$ and $2$,
whose hyperelliptic coverings have disjoint branch loci in $\P^1$. Let $X = X_1 \times_{\P^1} X_2$ and $Y$ the hyperelliptic curve ramified  over all branch points of $X_1$ and $X_2$.\\
Then the natural map $X \ra Y$ is an \'etale double covering whose Prym variety is isomorphic to $JX_1 \times JX_2$ as principally polarized abelian varieties.
\end{cor}
\end{ex}

\section{The Main Example}

Given two Prym-Tyurin varieties of exponent $q$ presented with
respect to the same group $G$, subgroup $H$ and rational
irreducible representations, but with group actions with disjoint
branch loci, Theorem \ref{prod} gives a new Prym-Tyurin variety of
exponent $[G:H] \cdot q$. Hence every Prym-Tyurin constructed with
a presentation with respect to a group action (and these are
almost all such varieties known up to now, see \cite{clrr}) gives
a new one. We will only give one example, which we think is
interesting,
since it arises also from a completely different geometric construction, as we shall see in the next section.\\

Let $p$ be an odd prime and
$$
G = \mathbf{D}_p = \langle \sigma, \tau : \sigma^p, \tau^2, (\sigma\, \tau)^2 \rangle
$$
the dihedral group of order $2p$.
Any complex
irreducible representation $V$ of degree two of $\mathbf{D}_p$ is defined over
the field $K_V = \mathbb{Q}(\zeta+\zeta^{-1})$, where $\zeta$ denotes a $p$-th root
of unity. The Galois group Gal$(K_V/\mathbb{Q})$ is cyclic of order
$\dfrac{p-1}{2}$, and the associated rational irreducible
representation $\mathcal{W}$ is of degree $p-1$ and given by
$$
\mathcal{W}: \; \sigma \rightarrow \left[ \begin{array}{ccccc}
0&& \cdots&0&-1\\
1&0&\cdots  &0&-1\\
0&1&\ldots & 0&-1\\
0 & 0 & \ddots\\
0&0 & \ldots  &1&-1\end{array} \right] \quad ,\quad \tau \rightarrow
\left[
\begin{array}{rrrrr}
0&0&\cdots &0&-1 \\
0&0&\cdots   &-1&0\\
0&0&\cdots &&0\\
0 & -1 & 0& \cdots & 0\\
-1&0&\cdots & & 0\end{array} \right]
$$
Apart from the trivial subgroup, $\mathbf{D}_p$ admits exactly 2
conjugacy classes of cyclic subgroups, namely the classes $C_1$ of
subgroups of order 2 and $C_2$ of subgroups of order $p$. Since
$\mathbf{D}_p$ can be generated by any even number $s \geq 4$ of
involutions with product equal to 1, there exist curves $Z$ with
$\mathbf{D}_p$-action and geometric signature $[0;(C_1,s)]$.

\begin{prop} \label{prop4.1}
Every $\mathbf{D}_p$-action on a curve $Z$ of geometric signature
$[0;(C_1,s)]$,  any subgroup $H$ of order 2 of $\mathbf{D}_p$ ,
and the representation $\cW$, give a presentation of the Jacobian
$JX$ (with $X = Z/H$) as a Prym-Tyurin variety of exponent $1$.
\end{prop}


\begin{proof}
Without loss of generality we may assume that $H = \{ 1 _G ,
\tau\}$. It suffices to show that the assumptions of Theorem
\ref{thm4.8} are satisfied with $r=1$ and $q=1$. Let the notation
be as in Section 2. Set $\pi_H:Z\to X$ and let $V$ denote a
complex irreducible representation associated to $\cW$ ($\deg V =
2$).

First we claim that $b=p$. To see this, note that any double coset
$HgH$ contains exactly $[H:H\cap g^{-1}Hg]$ right cosets of $H$.
Thus $H1_GH$  consists of one right coset, and for any $g$ not in
$H$ we have $[H:H\cap g^{-1}Hg]=|H|=2$. Hence the number $d =
|H\bs G/H|$ of double cosets is given by $p = [G:H]=1+2(d-1)$,
from where
$$
d=\frac{p+1}{2}.
$$
Moreover, $H\sigma^iH=H\sigma^jH$ if and only if $\sigma^i\in
H\sigma^jH=\{\sigma^j,\tau\sigma^j,\sigma^j\tau,\tau\sigma^j\tau=\sigma^{-j}\}$.
We conclude that $H\sigma^iH=H\sigma^jH$ if and only if $i =j$ or
$i=p-j$. This implies that a set of double coset representatives
for $H$ in $G$ is given by
$$
\left\{ g_{11}=1_G , g_{i1}=\sigma^{i-1} \; | \;  i=2, \ldots
,\frac{p+1}{2} \right\}.
$$
Therefore,
$$
b_1=\sum_{h \in H}
\chi_{\cW}(hg_{11}^{-1})=\chi_{\cW}(1_G)+\chi_{\cW}(\tau)=p-1
$$
and for $i=2, \ldots, \frac{p+1}{2}$
$$
b_i=\sum_{h \in H} \chi_{\cW}(hg_{i1}^{-1})=\chi_{\cW}(1_G\cdot
\sigma^{p-i+1})+\chi_{\cW}(\tau\cdot \sigma^{p -i+1})=-1,
$$
since $\chi_{\cW}(\tau\cdot \sigma^{i})=0$ and $\chi_{\cW}(\sigma^{i})=-1$, for all $i=1, \ldots ,p-1$.
This implies
$$
b= \gcd \left\{ b_1 - b_i \; | \; 2 \le i \leq \frac{p+1}{2}
\right\} = p.
$$
We conclude
$$
q= \frac{|\mathbf{D}_p|}{b \cdot \dim V} = 1.
$$

Consider $G_1=H$ as a representative of the conjugacy class $C_1$,
then we have that $\vert H\bs G/G_1 \vert =\dfrac{p+1}{2}$.

Due to Frobenius reciprocity,
$$
\dim V^{G_1}=\langle \text{Res}_{G_1} V,1_{G_1} \rangle
=\frac{1}{\dim V}\sum_{h\in G_1}\chi_V(h)=\frac{1}{2}(2+0)=1,
$$
and $G_1$ is a maximal subgroup with this property. It remains to
show that equation \eqref{eqn2.2} is satisfied for $q=1$. In fact,
$$
q[L_V:\Q](\dim V-\dim V^{G_1})-([G:H]-\vert H\bs G/{G_1}\vert
)=\frac{p-1}{2}(2-1)- \left( p-\frac{p+1}{2}\right)=0.
$$
The assertion now follows from Theorem \ref{prod}.
\end{proof}


As an immediate consequence of Proposition \ref{prop4.1}, Corollary \ref{cor3.3} and Theorem \ref{thm3.4} we obtain
\begin{prop} \label{prop4.2}
Set $G=\mathbf{D}_p$. For $i = 1$ and $2$ let $Z_i$ be a curve
with $G$-action of geometric signature $[0,(C_1,s_i)]$, with $s_i$
even integers $\geq 4$. Let $X_i$ be as in Proposition
\ref{prop4.1}. Assume moreover that the coverings $Z_i \ra \P^1$
have disjoint branch loci. Denote $Z = Z_1 \times_{\P^1} Z_2$ and
$X = X_1 \times_{\P^1} X_2$.

Then the action of the group $G^2$ on the curve $Z$ defines a
Prym-Tyurin variety $P$ in the Jacobian $JX$ of exponent $p$ and
dimension $\dfrac{p-1}{4}(s_1 + s_2 -8)$.

Moreover,
$$
P \simeq  JX_1 \times JX_2 \, .
$$
\end{prop}

 \medskip

In the sequel we use the following presentation of
$\mathbf{D}_p^2$:
$$
\mathbf{D}_p^2 = \; \langle \sigma_1, \tau_1, \sigma_2,\tau_2
\;|\; \sigma_i^p,\tau_i^2,(\sigma_i\tau_i)^2, [\sigma_1,\sigma_2],
[\tau_1,\tau_2], [\sigma_i,\tau_j]\; \forall\; 1 \leq i \neq j
\leq 2 \rangle .
$$
So $\sigma_1, \tau_1$ generate the first factor and $\sigma_2,
\tau_2$ generate the second factor of $\mathbf{D}_p^2$. Then, with
$H^2 = \langle \tau_1,\tau_2 \rangle$, we have
$$
Z_1 = Z/\langle \sigma_2,\tau_2 \rangle , \qquad Z_2 = Z/\langle
\sigma_1,\tau_1 \rangle ,
$$
$$
X_1 = Z/\langle \tau_1,\sigma_2,\tau_2 \rangle , \quad X_2 = Z/
\langle \sigma_1,\tau_1,\tau_2 \rangle \quad \text{and} \quad X =
Z/H^2.
$$
Moreover, consider the following quotient curves of $Z$:
\begin{itemize}
 \item $\widetilde{X} := Z/\langle \tau_1\tau_2\rangle$;
\item $\widetilde{Y}_j := Z/\langle
(\sigma_1^j\sigma_2,\tau_1\tau_2)\rangle$ for $j = 1, \ldots
\frac{p-1}{2}$; \item $Y := Z/\langle
\sigma_1,\sigma_2,\tau_1\tau_2\rangle$.
\end{itemize}
Then we have the following diagram with the degrees of the maps as indicated:
\begin{equation} \label{diag4.1}
\xymatrix@R=0.5cm{
& Z \ar[ddl]_{2p:1} \ar[d]^{2:1} \ar[ddr]^{2p:1} \\
& \widetilde{X} \ar[d]^{2:1} \ar[drr]^{p:1} \ar[drrrr]^{p:1} & & \\
Z_1 \ar[d]_{2:1} & X \ar[dl]^{p:1}_{q_1} \ar[dr]_{p:1}^{q_2} &
Z_2 \ar[d]^{2:1} & \widetilde{Y}_1 \ar[d]^{p:1}_{f_1}  & \dots & \widetilde{Y}_{\frac{p-1}{2}} \ar[dll]^{f_{\frac{p-1}{2}}}_{p:1}  \\
X_1 \ar[dr]^{p:1} && X_2 \ar[dl]_{p:1} & Y \ar[dll]^{2:1}  \\
& \P^1
}
\end{equation}

The next lemma gives the genus of the curves and the ramification
of the maps of this diagram, which we need in the sequel.

\begin{lem} \label{lem4.3}
{\em (a):} The covering $Z\to Y$ is \'etale, and hence so are
$Z\to \widetilde{X}$, $\widetilde{X} \to\widetilde{Y_j}$ and
$\widetilde{Y_j}\to Y$ for all $j$. $Y$ is hyperelliptic, ramified
over all $s_1 + s_2$ branch points of $Z \ra \P^1$. Hence
$$
g_Y=\frac{s_1+s_2}{2}-1 \quad \mbox{and} \quad g_{\widetilde{Y}_j}=\frac{p}{2}(s_1+s_2)-2p+1 \quad \mbox{for all} \; j;
$$

\medskip

{\em (b):} The covering $X_i \ra \P^1$, for $i = 1$ and $2$,  is
ramified exactly over the $s_i$ branch points of $Z_i \ra \P^1$,
each one with fibre consisting of $\frac{p-1}{2}$ points of order
$2$ and one unramified point. Hence
$$
g_{X_i}=\frac{s_1(p-1)}{4}-p+1 \quad \mbox{for} \; i = 1,2 \quad \text{ and } \quad g_X=\frac{(s_1+s_2)(p^2-p)}{4} -p^2 + 1.
$$
\end{lem}

\begin{proof}
The stabilizers of branch points of $Z \ra \P^1$ are conjugate
subgroups of either $\langle\tau_1\rangle$ or
$\langle\tau_2\rangle$, hence are of the form
$\langle\sigma_1^i\tau_1\rangle$ or
$\langle\sigma_2^i\tau_2\rangle$ for $i=0, \ldots, p-1$. For any
subgroup $U$ of $G$, the stabilizers for the covering $Z\to Z/U$
are given by the intersection of $U$ with the stabilizers for the
action of $G$. In our case the stabilizers for the action of
$\mathbf{D}_p^2$ are either trivial or are of order $2$, hence one
has to check whether the subgroup $U$ contains subgroups of the
form $\langle\sigma_1^i\tau_1\rangle$ or
$\langle\sigma_2^i\tau_2\rangle$. Doing this for the subgroups in
question, together with the Riemann-Hurwitz formula, gives the
assertions.
\end{proof}

\begin{lem} \label{lem4.4}
Let $Z \ra \P^1$ be a Galois covering with Galois group
$\mathbf{D}_p^2$ and geometric signature
$[0;(C_1^1,s_1),(C_1^2,s_2)]$ as above. Then

\noindent {\em (a):} $Y$ is the only intermediate hyperelliptic
curve of $Z \ra \P^1$ of genus $\frac{s_1+s_2}{2}-1$;

\noindent {\em (b):}  Up to isomorphism, the curves
$\widetilde{Y}_j$, $j=1, \ldots, \frac{p-1}{2}$ are the only
intermediate curves of $Z \ra \P^1$ which are \'etale $p-$fold
coverings of $Y$. They are pairwise non-conjugate in $Z \ra \P^1$;

\noindent {\em (c):} Up to isomorphism, $X$ is the only
intermediate curve of $Z \ra \P^1$ which is a $4:1$-quotient of
$Z$.
\end{lem}
\begin{proof}
This is a consequence of the subgroup structure of the group $\mathbf{D}_p^2$.\\
(a) follows from the fact that the subgroup $\langle\sigma_1,\sigma_2,\tau_1\tau_2\rangle$ is the only subgroup of index 2 of $\mathbf{D}_p^2$ not containing any stabilizer for the action of $\mathbf{D}_p^2$ on $Z$.\\
(b) First note that the subgroup
$\langle\sigma_1^j\sigma_2,\tau_1\tau_2\rangle$ corresponding to
$\widetilde{Y}_j$ is a subgroup of index $p$ in the subgroup
corresponding to $Y$. According to Lemma \ref{lem4.3} (a),
$\widetilde{Y}_j \ra Y$ is \'etale of degree $p$. Moreover, note
that whenever we have a curve $Z$ with group action by a group
$G$, the covering $Z\to Z/G$ is the Galois closure of an
intermediate covering $Z/U\to Z/G $ if and only if the core
$\Core_G U:=\bigcap_{g \in G} U^g$ is trivial. In our situation,
the only subgroups, up to conjugacy, of $\mathbf{D}_p^2$ of index
$2p$ with trivial core in $\mathbf{D}_p^2$
are the ones determining the $\widetilde{Y_j}$. It is easy to see that they are pairwise non-conjugate.\\
(c): In $\mathbf{D}_p^2$ there is just one conjugacy class of
subgroups of order $4$. All of them are non-cyclic. Therefore the subgroup defining
$X$ is uniquely determined up to conjugacy in $\mathbf{D}_p^2$.
\end{proof}
As an immediate consequence we get
\begin{cor} \label{cor4.5}
Suppose $\widetilde{Y} \ra Y$ is an \'etale $p$-fold covering of a
hyperelliptic curve $Y$ such that the Galois closure $Z$ of the
composed map $\widetilde{Y} \ra Y \ra \P^1$ has Galois group
$\mathbf{D}_p^2$. Then $\widetilde{Y}$ is (isomorphic to) one of
the curves $\widetilde{Y}_j$ of diagram \ref{diag4.1}. In
particular it determines the curves $Z_1$ and $Z_2$ of that
diagram uniquely and $X, X_1$ and $X_2$ uniquely up to conjugacy.
 \end{cor}


\section{\'Etale $p$-fold coverings of hyperelliptic curves}

In this section we give a geometric construction of the $p$-fold coverings $X_1$ and $X_2$ of diagram \ref{diag4.1} in terms of the covering $f: \widetilde{Y} \ra Y$.

\subsection{The curves}
So let $f: \widetilde{Y} \ra Y$ denote an \'etale covering of degree $p$ of a hyperelliptic curve $Y$. Suppose the (2:1)-covering $h: Y \ra \P^1$ is ramified over $B_h = \{a_1, \ldots, a_{2g+2}\} \subset \P^1$.

Define the curve $X$ as the following fibre product,
\begin{equation} \label{diag5.1}
\xymatrix{
  X:=(f^{(2)})^{-1}(g_2^1) \ar[d]_{\pi=f^{(2)}|_X} \ar@{^{(}->}[r]
                & \widetilde{Y}^{(2)} \ar[d]^{f^{(2)}}  \\
  \P^1\cong g_2^1  \ar@{^{(}->}[r]   & Y^{(2)}.             }
\end{equation}
Here $Y^{(2)}$ and $\widetilde{Y}^{(2)}$ denote the
second symmetric product of $Y$ and $\widetilde{Y}$.
In order to work out conditions for $X$ to be smooth and irreducible, we need some notation.

Fix a point $z_0 \in \P^1 \setminus B_h$ and let $\gamma_i$ denote
the class of a path at $z_0$ going once around $a_i$ clockwise as
usual, then the fundamental group of $\P^1 \setminus B_h$ is
\begin{equation} \label{eqn5.2}
\pi_1(\P^1 \setminus B_h, z_0)=\;\langle \gamma_1, \ldots
,\gamma_{2g+2}: \prod_1^{2g+2}\gamma_i=1 \rangle
\end{equation}
If $\iota$ is the hyperelliptic involution of $Y$, we denote
$$
h^{-1}(z_0)=\{x,\iota x\},
$$
$$
f^{-1}(x)=\{x_1, \dots, x_p\}\subset \widetilde{Y} \quad
\mbox{and} \quad f^{-1}(\iota x)=\{y_1, \dots, y_p\}\subset
\tilde{Y}.
$$
If $\mu:\pi_1(\P^1 \setminus B_h, z_0) \to S_{2p}$ is a classifying
homomorphism of  $h\circ f:\tilde{Y} \to Y \to \P^1$, we know that
$$
\mu(\gamma_i)=t_1t_2\dots t_p
$$
where $t_1,t_2,\dots ,t_p$ are disjoint transpositions. But not all such products can occur.
In fact, if we identify for $i = 1, \ldots,p$,
\begin{equation}\label{eqn5.3}
 x_i = i \qquad \mbox{and} \qquad y_i = p+i,
\end{equation}
then we have the following lemma, the proof of which is obvious.
\begin{lem}
A homomorphism $\mu:\pi_1(\P^1 \setminus B_h, z_0) \to S_{2p}$ is
a classifying homomorphism for  $h\circ f:\widetilde{Y} \to Y \to
\P^1$ if and only if
\begin{enumerate}
 \item $G := \Ima(\mu)$ is an imprimitive transitive subgroup of imprimitivity degree $p$ of $S_{2p}$ and
\item $\mu(\gamma_i) = t_1 t_2 \cdots t_p$ with disjoint
transpositions $t_i$ of the form $(j \;\;p+k)$ with $1 \leq j,k
\leq p$.
\end{enumerate}
\end{lem}
Finally, if we denote for $i,j = 1, \ldots, p$,
\begin{equation} \label{eqn5.4}
P_{ij} := x_i + y_j \in X \subset \widetilde{Y}^{(2)},
\end{equation}
then $\pi^{-1}(z_0) = \{ P_{ij} \; | \; i,j = 1, \ldots, p \}$,
with $\pi : X \to \mathbb{P}^1$ from Diagram \ref{diag5.1}.\\

Now consider the group
$\G \subset S_{2p}$
generated by
\begin{itemize}
\item $\varphi_1=\Pi_{i=1}^p(i\;\; i+p)$
\item $\varphi_2=[\Pi_{i=1}^{p-1}(i\;\; i+p+1)](p\;\; p+1)$
\item $\varphi_3=(1\;\; p+1)\Pi_{i=2}^p(i\;\; 2p+2-i)$ and
\item $\varphi_4=(1\;\; p+2)(2\;\; p+1)\Pi_{i=3}^p(i\;\; 2p+3-i)$
\end{itemize}

\begin{lem}  \label{lem5.1}
\begin{equation} \label{eqn5.5}
\G = \mathbf{D}_p \times \mathbf{D}_p \subset S_{2p},
\end{equation}
where $\varphi_1,\varphi_2$ generate the first factor $\mathbf{D}_p$ and
$\varphi_3,\varphi_4$ the second.
\end{lem}

\begin{proof} Obviously all $\varphi_i$'s are of order $2$. Moreover
we have the following relations:
$$
\varphi_1\cdot \varphi_2 = (1\;\; p\;\; p-1\;\; p-2\; \dots \; 2)(p+1\;\;p+2\; \dots \; 2p),
$$
$$
\varphi_3\cdot \varphi_4 =(1\; 2\; 3\; 4\; \dots \; p)(p+1\;p+2\; \dots \; 2p).
$$
Hence $\langle \varphi_1,\varphi_2 \rangle$ and $\langle
\varphi_3,\varphi_4 \rangle$ are both groups isomorphic to
$\mathbf{D}_p$. Moreover, one easily checks that $\varphi_1$ and
$\varphi_2$ commute with $\varphi_3$ and $\varphi_4$, which
completes the proof.
\end{proof}

\begin{prop} \label{prop5.2}
If $f: \widetilde{Y} \to Y$ is an \'etale covering of degree $p$
of a hyperelliptic curve $Y$ such that the image of a classifying
homomorphism $\mu: \pi_1(\P^1 \setminus B_h, z_0) \to S_{2p}$ is
the group $\G$, then the curve $X$ of diagram \ref{diag5.1} is
smooth and irreducible.
\end{prop}

\begin{proof}
The stabilizer of the element $P_{11} = x_1 + y_1$ of the fibre
$\pi^{-1}(z_0)$ of the map $\pi: X \ra \P^1$ is the group
$$\G_{P_{11}} =\; \langle \varphi_1,\varphi_3 \rangle,
$$
which is Klein's group of 4 elements. Since $\G$ is of order $4p^2$,
this means that $\G$ acts transitively on the set $\{P_{ij} \;|\; i,j = 1, \ldots, p\}$ implying that $X$ is irreducible. The proof
of the fact that $X$ is smooth is a slight
generalization of the proof of \cite{lb}, Lemma 12.8.1. (see also \cite[Lemma 3.1]{LRR}, where the special case $p=3$ is proved).
\end{proof}

\begin{prop} \label{prop5.3}
The curve $X$ of diagram \eqref{diag5.1} coincides with the curve $X$ in diagram \eqref{diag4.1}.
\end{prop}
\begin{proof}
According the the proof of Proposition \ref{prop5.2}, the curve
$X$ corresponds to the Klein subgroup $\langle \varphi_1,
\varphi_3\rangle$ of $\mathbf{D}_p^2$. In Lemma \ref{lem4.4} we
show that, up to conjugacy, there is only one such a subgroup.
\end{proof}

\subsection{The correspondence}

In the last subsection we saw how to describe the curve $X$ of
diagram \ref{diag4.1} in terms of the covering $f : \widetilde{Y}
\ra Y$. We can use this to define a fixed-point free symmetric
effective $(p-1)$-correspondence $D$ on $X$.

For this fix a point $z_0 \in U := \P^1 \setminus B_h$ and denote
the fibre $\pi^{-1}(z_0) = \{ P_{ij} \;|\; i,j = 1, \ldots p \}$
as in equation \eqref{eqn5.4}. Moreover we use the notation
$$
I_{ij} := \{ (k,l) \in \{1, \ldots , p\}^2 \; | \; k+l \not \equiv i +j \mod p \quad \mbox{and} \quad k-l \not \equiv i-j \mod p \}.
$$
Then in the fibre $\pi^{-1}(z_0)$ the correspondence $D$ is defined by
\begin{equation}  \label{eqn5.6}
D_{z_0}(P_{ij}) := \sum_{(k,l) \in I_{ij}} P_{kl}
\end{equation}
We extend $D_{z_0}$ to a correspondence $D_U$ on $\pi^{-1}(U)$ in
the usual way (see e.g. \cite[Section 3]{LRR}) as follows: We
enumerate the $x_i$ and $y_j$ in such a way that the stabilizer of
$P_{11}$ is the group $H^2 =\; \langle \varphi_1,\varphi_3
\rangle$. If $\{ (g_{ij},g_{kl}) \;|\; 1 \leq i,k \leq d, 1 \leq
j,l \leq n_i \} $ denotes the set of representatives of the right
and left cosets of $H^2$ in $\G$ as in the proof of Lemma
\ref{lem3.6} with $d = \frac{p+1}{2}, n_1 = 1$ and $n_i = 2$ for
$i \geq 2$ (see proof of Proposition \ref{prop4.1}), this induces
a $\G$-equivariant bijection
$$
\{ (g_{11},g_{11}), \ldots, (g_{d n_d}, g_{d n_d}) \}
\longrightarrow \pi^{-1}(z_0) = \{ P_{ij} \;|\; i,j = 1, \ldots p
\}
$$
to be described in the proof of Proposition \ref{prop6.1}.

Then for every point $z \in U$ choose a path $\gamma_{z}$ in $U$
connecting $z$ and $z_0$. The path defines a bijection
$$
\mu: \pi^{-1}(z) \to \pi^{-1}(z_0)=\{P_{11}, \dots, P_{pp}\}
$$
in the following way: For any $x \in \pi^{-1}(z)$ denote by
$\tilde{\gamma}_x$ the lift of $\gamma_{z}$ starting at $x$. If
$P_{ij} \in \pi^{-1}(z_0)$ denotes the end point of
$\tilde{\gamma}_x$, define
$$
\mu(x) = P_{ij}.
$$
Define
$$
D_U := \{ (x, x') \in \pi^{-1}(U) \times_{\P^1} \pi^{-1}(U) \;|\;
\mu(x') \in D_{z_0}(\mu(x)) \}.
$$
Finally, define $D$ to be the closure of $D_U$ in $X\times X$.

\begin{prop}\label{kanevD}
$D$ is a correspondence on $X$.
\end{prop}

In the next section we will see that $D$ coincides with the Kanev
correspondence associated by the correspondence ${\cD}$ defining
the Prym-Tyurin variety $P$ of Proposition \ref{prop4.2}. In
particular $D$ is an  effective symmetric fixed-point free
correspondence whose associated endomorphism $\gamma_D$ on the
Jacobian $JX$ satisfies the equation  $\gamma_D^2 + (p-2)\gamma_D
-(p-1) = 0$. Of course it is easy to see this also directly.

\begin{proof}
We have to show that the definition of $D_U$ is independent of the
choice of the paths $\gamma_z$. For this it suffices to show that
$D_U$ is invariant under the action of $\G=\mathbf{D}_p\times
\mathbf{D}_p$; that is, the diagram

$$
\xymatrix{
 P_{ij} \ar[d]_D \ar[rr]^{\varphi_k}
               & & \varphi_k(P_{ij}) \ar[d]^D  \\
 D(P_{ij})  \ar[rr]_{\varphi_k}  & & D(\varphi_k(P_{ij}))=\varphi_k(D(P_{ij}))             }
$$
commutes for $1 \leq k \leq 4$. In fact, for $\varphi_1$ we have
$$
\varphi_1^{-1}(D(\varphi_1(P_{ij}))) = \varphi_1(D(P_{ji}))  =
\varphi_1 \left( \sum_{(l,k) \in I_{ji}} P_{lk}\right) =
\sum_{(k,l) \in I_{ij}} P_{kl} = D(P_{ij}).
$$
For the other $\varphi_k$ the proof is similar, using
$$
\varphi_2(P_{ij})=P_{kl} \quad \mbox{with} \quad k \equiv j-1\mod p, \quad l \equiv i+1\mod p,
$$
$$
\varphi_3(P_{ij})=P_{kl} \quad \mbox{with} \quad k \equiv p-j+2 \mod p, \quad  l \equiv p-i+2 \mod p,
$$
$$
\varphi_4(P_{ij})=P_{kl} \quad \mbox{with} \quad k \equiv p-j+3 \mod p, \quad  l \equiv p-i+3 \mod p.
$$\end{proof}

\begin{ex}
In the special case $p=3$ the curve $X$ coincides with the corresponding curve in \cite{LRR}. Moreover we have
$$
D(P_{ij}) = \sum_{(k,l) \in I_{ij}} P_{kl} = \sum_{k=1, k \neq i}^3 P_{kj} + \sum_{l=1, l \neq i}^3 P_{il}.
$$

This is just the correspondence of \cite{LRR}. Hence our
construction of Prym-Tyurin varieties is a generalisation to
arbitrary odd prime $p$ of the  Prym-Tyurin varieties of
\cite{LRR} for $p=3$.
\end{ex}

\section{Comparison of the correspondences ${\cD}$ and $D$}

According to Proposition \ref{prop5.3} the curve $X$ of diagram
\eqref{diag5.1} coincides with the curve $X$ of diagram
\eqref{diag4.1}.
Hence we have two correspondences on $X$, namely ${\cD}$ defined
by equation \eqref{eqn3.6} (in the special case $G = \mathbf{D}_p$)
and $D$ defined in the previous section. In this section
we compare both and show that they induce the same Prym-Tyurin variety in $JX$.\\

Recall from   \cite{clrr} that for a general ${\cD}$ as in
equation \eqref{eqnD}, we constructed  a correspondence
$\cK_{{\cD}}$, called the associated {\it Kanev correspondence},
which is effective, symmetric, fixed-point free and whose
associated endomorphism $\gamma_{\cK_{{\cD}}}$ satisfies the
equation $\gamma_{\cK_{{\cD}}}^2 + ({q} -2) \gamma_{\cK_{{\cD}}} -
({q} -1) = 0$.  In fact, in the proof of \cite[Lemma 3.8]{clrr} we
saw that $\cK_{{\cD}}(x)$ is given by
\begin{equation} \label{eqn6.1}
\cK_{{\cD}}(x)
=\sum_{i=2}^d\left[\frac{b_1-b_i}{b}-1\right]\sum_{j=1}^{n_i} \pi_H g_{ij}(z)
\end{equation}
for all $x \in X$ and $z \in Z$ with $\pi_{H}(z) = x$, where
$$
b_i := \sum_{k=1}^r \sum_{h \in H} \tr_{L/\QQ}(\chi_{V_k}(hg_{i1}^{-1}))
$$
Moreover, in terms of $\cK_{{\cD}}$ the associated Prym-Tyurin variety $P$ is given by
$$
P = \Ima(1-\gamma_{\cK_{{\cD}}}).
$$

Now let the notation be as in Section 5. In particular $Z$ is a
Galois covering of $\P^1$ with Galois group $\mathbf{D}_p^2$ and
$\pi_{H^2}: Z \ra X$ is the covering of diagram \eqref{diag5.1}.
Let $\cD$ denote the correspondence defined in equation
\eqref{eqn3.6} in our special situation. Then we have
\begin{prop} \label{prop6.1}
The Kanev correspondence $\cK_{\cD}$ associated to $\cD$
coincides  with the correspondence $D$ of Proposition
\ref{kanevD}.
\end{prop}
\begin{proof}
For the proof  we show that the right hand side of \eqref{eqn6.1}
in this special situation is equal to the right hand side of
\eqref{eqn5.6} under the identifications given below.

For this we recall the notations. The underlying group is
$G^2=\mathbf{D}_p\times \mathbf{D}_p$ with generators of the two
factors $\sigma_i$ of order $p$ and $\tau_i$ of order 2 for $i=1$
and $2$. We are considering the representations $\cW\otimes V_0$
and  $V_0\otimes \cW$, with $\cW$ as in Proposition \ref{prop4.1}.
Therefore here $r=2$.

We need to collect some of the previous results for $G =
\mathbf{D}_p$ and its subgroup $H=\langle \tau \rangle$ from
Proposition \ref{prop4.1}:
\begin{itemize}
\item  the double cosets representatives for $H\bs G/H$ we
consider are $\sigma^{i-1}$, with $i=1, \ldots ,\frac{p+1}{2}$;

\item the coefficients of the correspondence for $\mathbf{D}_p$
and $H$ on $X_i$ are $a_1=p-1$ and $a_i=-1$ for $i=2, \ldots
,\frac{p+1}{2}$;

\item the double coset represented by $1_G$ has $2$ elements,
hence just one right coset. The double coset represented by
$\sigma^i$ for $i > 1$  has $4$ elements, hence two right cosets.
The representatives for the right and left cosets in this double
coset are $\sigma^i$ and $\sigma^{p-i}$.
\end{itemize}

We now work out the induced representatives and coefficients for
the correspondence for $G^2=\mathbf{D}_p\times \mathbf{D}_p$ and
its subgroup $H^2=\langle \tau_1 \rangle \times \langle
\tau_2\rangle$ on the curve $X$.

\begin{itemize}

\item the double coset representatives for $H^2\bs G^2/H^2$ we
consider are $(\sigma_1^{i-1},\sigma_2^{j-1})$  for $i,j=1,
\ldots,\frac{p+1}{2}$;

\item the representatives for the right and left cosets inside the
double coset are:

for $i=j=1$ there is only 1 right coset represented by
$(1_G,1_G)$;

for $i=1$ and $1 \leq k \leq \frac{p-1}{2}$ the right and left
representatives are $\{(1_G,\sigma_2^{k}),(1_G,\sigma_2^{p-k})\}$.

for $1 \leq \ell \leq \frac{p-1}{2}$ and $j = 1$ they are
$\{(\sigma_1^{\ell},1_G), (\sigma_1^{p-\ell},1_G)\}$.

for $1 \leq \ell , k  \leq \frac{p-1}{2}$ they are
$\{(\sigma_1^{\ell},\sigma_2^k),(\sigma_1^{p-\ell},\sigma_2^k),
(\sigma_1^{\ell},\sigma_2^{p-k}),
(\sigma_1^{p-\ell},\sigma_2^{p-k})\}$.
\end{itemize}

Now consider the defining equation \eqref{eqn3.6} for the
correspondence $\cD$. According to equations \eqref{eqn3.8} and
\eqref{eqn3.9} its coefficients $b_{ij} = a_{ij} + a_{ij}'$ (where
the $a_{ij}$ correspond to $\cW \otimes V_0$ and the $a_{ij}'$
correspond to $V_0 \otimes \cW$) are given by
$$
a_{1j}=|H|a_1=2(p-1) \;  \text{for all} \; j;
$$
$$
a_{ij}=|H|a_i=-2 \; \text{ for all } i = 2, \ldots, \frac{p+1}{2}
\; \mbox{and all} \; j;
$$
$$a_{i1}'=|H|a_1=2(p-1) \;  \text{for all}  \; i;
$$
$$
a_{ij}'=|H|a_i=-2  \; \text{ for all }\;  i \; \mbox{and} \; j =
1, \ldots, \frac{p+1}{2}.
$$
Hence we get
$$
b_{ij} = \left\{ \begin{array}{cll}
                 4p-4 &  & i=j=1;\\
                 2p-4 & for & i=1,\; j=2, \cdots, \dfrac{p+1}{2}\; \text{ and }\;
                           i = 2 , \ldots, \frac{p+1}{2},\; j = 1;\\
                 -4 &  & 2 \leq i, j \leq \frac{p+1}{2}.
                \end{array} \right.
$$
and $b= \gcd\{b_{11} - b_{ij} \} = 2p$.

\medskip

Therefore, from equation \eqref{eqn6.1} we see that the
coefficients for the correspondence $\cK_{\cD}$ vanish if either
$i=1$ or $j=1$, and are equal to $1$ for $i,j\geq 2$. This gives
\begin{equation} \label{eqn6.2}
\cK_{\cD}(x)=\sum_{\ell,k=1}^{\frac{p-1}{2}}[\pi_{H^2}(\sigma_1^{\ell},\sigma_2^k)(z)
+\pi_{H^2}(\sigma_1^{\ell},\sigma_2^{p-k})(z)+\pi_{H^2}(\sigma_1^{p-{\ell}},\sigma_2^k)(z)
+\pi_{H^2}(\sigma_1^{p-{\ell}},\sigma_2^{p-k})(z)].
\end{equation}
for any $x\in X$ and $z\in Z$ mapping to $x$, and where
$\pi_{H^2}: Z \ra X$ is the covering corresponding to the subgroup
$H^2$.

Now the identification of the two definitions of the curve $X$ is
as follows. If $P_{11}$ corresponds to the point $x \in X$, then
$P_{\ell k}$ corresponds to the point
$\pi_{H^2}(\sigma_1^{\ell},\sigma_2^k)(z)$ where $z \in Z$ maps to
$x$. This follows from the fact that the identification is
$G^2$-equivariant. This means that for $P_{rs}\in X$ we have
$$
\pi_{H^2}(\sigma_1^{\ell},\sigma_2^k)(P_{rs})=P_{uv}
$$
with $u$ and $v$ such that $u+v \equiv r+s+{\ell} \mod p$ and
$u-v \equiv r-s+k \mod p$. Inserting this in \eqref{eqn6.1}, we
obtain
$$
\cK_{\cD}(P_{ij}) = \sum_{(k,l) \in I_{ij}} P_{kl},
$$
which completes the proof of the proposition.
\end{proof}

As an immediate consequence of Propositions \ref{prop4.2} and  \ref{prop6.1} we obtain

\begin{cor} \label{cor6.2}
Suppose $\widetilde{Y} \ra Y$ is an \'etale $p$-fold covering of a hyperelliptic curve $Y$ such that the Galois closure $Z$ of the composed map $\widetilde{Y} \ra Y \ra \P^1$ has Galois group $\G = \mathbf{D}_p^2$.\\
{\em (a):} The correspondence $D$ on the curve $X$ of diagram \eqref{diag5.1} defines a Prym-Tyurin variety $P$ of exponent $p$ in the Jacobian $JX$.\\
{\em (b):} There exist curves $X_1$ and $X_2$ over $\P^1$, whose fibre product over $\P^1$ is $X$ such that
$$
P \simeq JX_1 \times JX_2.
$$
as principally polarized abelian varieties.
\end{cor}

\section{Decomposition for the Jacobian of $X$}\label{S:desc}

Let the notation be as in Section 4. So we are given the curve $Z
= Z_1 \times_{\P^1} Z_2$ with $\mathbf{D}_p^2$-action as in
Proposition \ref{prop4.2}. Consider the curve $X$ of diagram
\eqref{diag4.1} (or equivalently the curve $X$ of diagram
\eqref{diag5.1}). If $X_i, Y$ and $\widetilde{Y}_j$ are as in
diagram \eqref{diag4.1}, recall that $\widetilde{Y_j} \ra Y$ are
\'etale $p$-fold coverings. If then $P(\widetilde{Y}_j/Y)$ denotes
the corresponding (generalized) Prym variety (i.e. the connected
component of the origin of the kernel of the norm map
$J\widetilde{Y}_j \ra JY$), the Jacobian $JX$ decomposes up to
isogeny in the following way.

\begin{thm} \label{thm7.1}
$$
JX \sim JX_1 \times JX_2 \times P(\widetilde{Y}_1/Y) \times \cdots \times P(\widetilde{Y}_{(p-1)/2}/Y).
$$
\end{thm}

For the proof we use some results of \cite{CR} which we recall
first: For a given finite group $G$, let $Z \ra \P^1$ denote a
Galois covering with Galois group $G$, $H$ a subgroup of $G$ and
$Z_{H} := Z/H$. To every rational irreducible representation $\cW$
of $G$ one can associate an abelian subvariety $B_{\cW}$ of the
Jacobian $JZ$ which is uniquely determined up to isogeny. Let
$V_0$ denote the trivial representation and $\cW_1, \ldots ,
\cW_s$ the non-trivial rational irreducible representations of
$\G$ with associated complex irreducible representations $V_1,
\ldots, V_s$. Moreover,  $\rho_{H}$ denotes the character of $G$
induced by the trivial character of $H$. Then \cite[Lemma 4.3]{CR}
says that
$$
\rho_{H} = V_0  \oplus \bigoplus_{j=1}^s c_j \mathcal{W}_j \quad
\mbox{with} \quad c_j = \frac{\dim V_j^H}{m_j},
$$
where $m_j$ denotes the Schur index of the representation $V_j$.
Moreover, \cite[Proposition 5.2]{CR} says that the Jacobian
$JZ_{H}$ decomposes up to isogeny as
\begin{equation} \label{eqn7.1}
JZ_{H} \sim B_{\cW_1}^{c_1} \times \cdots \times B_{\cW_s}^{c_s}.
\end{equation}
Finally, if $H \subset N$ are two subgroups of $G$, then
$$
\rho_{H} = \rho_{N} \oplus \bigoplus_{j=1}^s d_j \mathcal{W}_j
\quad \mbox{with} \quad d_j = \frac{\dim V_j^H - \dim V_j^N}{m_j}
\; ( \geq 0)
$$
and, according to \cite[Corollary 5.4]{CR},  the Prym variety
$P(Z_{H}/Z_{N})$ of the morphism $ Z_{H} \ra Z_{N}$ decomposes up
to isogeny as
\begin{equation} \label{eqn7.2}
P(Z_{H}/Z_{N}) \sim B_{\cW_1}^{d_1} \times \cdots \times
B_{\cW_s}^{d_s}.
\end{equation}
\vspace{0.3cm}

Returning to the dihedral group $\mathbf{D}_p$, we denote by
$V_j$, $j=1,\dots ,\frac{p-1}{2}$ its complex irreducible
representations of degree $2$, with corresponding characters given
by
$$
\chi_{V_j}(\sigma^h) = \omega^{jh} + \omega^{-jh}, \qquad
\chi_{V_j}(\tau\sigma^h) = 0,
$$
where $\omega$ denotes a fixed $p$-th root of unity.
Then $\cW=\bigoplus_{j=1}^{\frac{p-1}{2}}V_j$ is the irreducible rational representation of degree $p-1$.

Consider, as in Section 4, the group $\mathbf{D}_p^2 = \;
<\sigma_1, \tau_1,\sigma_2.\tau_2>$. The complex irreducible
representation of $\mathbf{D}_p^2$ are the (outer) tensor products
of the irreducible representations of the factors. We consider the
following rational irreducible representations of
$\mathbf{D}_p^2$: the trivial $\cV_0 = V_0 \otimes V_0$, the
alternating representation $\cV_0' = V_0' \otimes V_0'$,
$$
\cW_1 = \cW \otimes V_0 \quad \mbox{and} \quad \cW_2 = V_0\otimes
\cW ,
$$
and, for $j = 1, \ldots ,\frac{p-1}{2}$,
$$
\cU_j = \bigoplus_{i=1}^{\frac{p-1}{2}} (V_i \otimes V_k)= (V_1
\otimes V_j)\oplus (V_2 \otimes V_{j+1}) \oplus \ldots  \oplus
(V_{(p-1)/2} \otimes V_{j-1})
$$
where $1 \leq k \leq \frac{p-1}{2}$ is given by $j+i-1$ if $j+i-1
\leq (p-1)/2$ and $k = j+i-1-(p-1)/2$ otherwise.

Recall the subgroups defining the curves of diagram \ref{diag4.1}:

\begin{itemize}
\item $X$ is defined by $H^2:=\;\langle\tau_1,\tau_2\rangle$,
\item $X_1$ is defined by
$H_1:=\;\langle\sigma_2,\tau_1,\tau_2\rangle$, \item $X_2$ is
defined by $H_2:=\;\langle\sigma_1,\tau_1,\tau_2\rangle$, \item
$\widetilde{Y}_j$ is defined by
$L_j:=\;\langle\sigma_1^j\sigma_2,\tau_1\tau_2\rangle$, \item $Y$
is defined by $M:=\;\langle\sigma_1,\sigma_2,\tau_1\tau_2\rangle$,
\end{itemize}
\begin{lem} \label{lem7.2}
{\em (a):} $\rho_H = \cV_0 \oplus \cW_1 \oplus \cW_2 \oplus \cU_1 \oplus \cdots \oplus \cU_{\frac{p-1}{2}},$ \\
{\em (b):} $\rho_{H_i} = \cV_0 \oplus \cW_i$ for $i = 1$ and $2$,\\
{\em (c):} $\rho_M = \cV_0 \oplus \cV_0'$,\\
{\em (d):} $\rho_{L_j}=\rho_M\oplus \cU_j$ for $ j = 1, \ldots, \frac{p-1}{2}$.
\end{lem}

\begin{proof}
(a): As $H^2$ is a subgroup of $\mathbf{D}_p^2$ satisfying Hypothesis \ref{hyp},
we have for the character product
$$
\langle\rho_{H^2}, \cW_1 \rangle \;= \;\langle\rho_{H^2}, \cW_2
\rangle \;= 1.
$$
By direct computation of the character product   using Frobenius
reciprocity, one shows  $\langle \rho_{H^2}, \cU_j \rangle = 1 $
for all $j$. Finally, as the sum of the degrees of these
representations is the degree of $\rho_{H^2}$, we get the result.
The proofs of the other assertions are similar.
\end{proof}

\begin{proof} {\it of Theorem \ref{thm7.1}}.
Let $B_{\cW_i}$ and $B_{\cU_j}$ denote the abelian varieties associated to the representations $\cW_i$ and $\cU_j$. Then we have, according to Lemma \ref{lem7.2} (a) and equation \eqref{eqn7.1},
$$
JX \sim B_{\cW_1} \times B_{\cW_2} \times B_{\cU_1} \times \cdots \times B_{\cU_{\frac{p-1}{2}}.}
$$
It remains to identify the $B$'s. According to Lemma \ref{lem7.2} (b) and equation \eqref{eqn7.1},
$$
B_{\cW_i} \sim JX_i \quad \mbox{for} \quad i = 1,2
$$
and Lemma \ref{lem7.2} (c) and (d) and equation \eqref{eqn7.2} imply
$$
B_{\cU_j} \sim P(\widetilde{Y}_j/Y) \quad \mbox{for all} \quad j.
$$
This completes the proof of the theorem.
\end{proof}
\bibliographystyle{amsplain}

\end{document}